\documentclass[12pt, reqno]{amsart}
\usepackage{amsmath, amstext, amsbsy, amssymb, amscd}
\usepackage{amsmath}
\usepackage{amsxtra}
\usepackage{amscd}
\usepackage{amsthm}
\usepackage{amsfonts}
\usepackage{amssymb}
\usepackage{eucal}
\usepackage{color}

\newtheorem{theorem}{Theorem}[section]

\newtheorem{lemma}[theorem]{Lemma}
\newtheorem{proposition}[theorem]{Proposition}
\newtheorem{corollary}[theorem]{Corollary}
\theoremstyle{definition}
\newtheorem{definition}[theorem]{Definition}
\newtheorem{example}[theorem]{Example}

\theoremstyle{remark}
\newtheorem{remark}[theorem]{Remark}

\definecolor{A}{rgb}{.75,1,.75}

\numberwithin{equation}{section}

\newcommand{\An}{\dahc}
\newcommand{\Bn}{\sdaha}
\newcommand{\Anhat}{\mathfrak{H}^{\mathfrak c}}
\newcommand{\Bnhat}{ \mathfrak{H}^-}

\newcommand{\C}{ \mathbb C }
\newcommand{\Cl}{ {\mathcal C} }
\newcommand{\ind}{\text{Ind}}

\newcommand{\N}{\mathbb N}
\newcommand{\tS}{\widetilde{S}_}
\newcommand{\tSn}{\widetilde{S}_n}
\newcommand{\Tn}{\C S_n^-}
\newcommand{\Z}{ \mathbb Z }
\newcommand{\Zn}{ \mathcal Z}

\newcommand{\dahc}{\ddot{\mathfrak H}^{\mathfrak c}}             
\newcommand{\sdaha}{\ddot{\mathfrak H}^-}             

{\vskip-\lastskip\medskip
  \noindent
  {\em #1.}\enspace
  }%
{\qed\par\medskip
  }

\begin{document}
\title[Spin double affine Hecke algebras]
{Double affine Hecke algebras for the spin symmetric group}

\author[Weiqiang Wang]{Weiqiang Wang}
\address{Department of Mathematics, University of Virginia,
Charlottesville, VA 22904} \email{ww9c@virginia.edu}


\begin{abstract}
We introduce a new class (in two versions, $\An$ and $\Bn$) of
rational double affine Hecke algebras (DaHa) associated to the
spin symmetric group. We establish the basic properties of the
algebras, such as PBW and Dunkl representation, and connections to
Nazarov's degenerate affine Hecke-Clifford algebra and to a new
degenerate affine Hecke algebra introduced here. We formulate a
Morita equivalence between the two versions of rational DaHa's.
The trigonometric generalization of the above constructions is
also formulated and its relation to the rational counterpart is
established.
\end{abstract}

\maketitle

\date{3/27/06}


\section{Introduction}
\subsection{}

The double affine Hecke algebras (DaHa) introduced by Cherednik \cite{Ch}
are intimately related to K.~Saito's elliptic root systems \cite{Sai}
and the Dunkl operators \cite{Dun}, and they have
numerous connections and applications to
Macdonald polynomials, integrable systems, and (quantum) affine algebras,
etc. The DaHa affords three variants of degeneration:
rational, trigonometric, and elliptic. The rational DaHa, which
goes back as a special case to Drinfeld \cite{Dr},
has been actively studied in recent years by
many authors (see Etingof-Ginzburg \cite{EG}
and the review of Rouquier \cite{Rou} for extensive references).

The question addressed here is whether or not a reasonable notion
of DaHa's associated to the spin symmetric group of Schur
\cite{Sch} exists. For Coxeter groups, there is a standard
procedure to construct the associated Hecke algebras. The spin
symmetric group is not a reflection group, and so a priori it is
not clear whether such a DaHa should exist and, if it exists, how
its characteristic feature looks like.

In this paper we provide a natural construction of rational and
trigonometric DaHa's, denoted by $\An$ and $\Bn$, associated to
the spin symmetric group, which not only exhibits new phenomenon
but also suggests for further natural generalization. These new
algebras afford favorable properties similar to the usual DaHa
(see \cite{EG} and Suzuki \cite{Suz}). Instead of constructing
directly the DaHa for the spin symmetric group, our idea is to
start with a construction of DaHa for a certain double cover of
the hyperoctahedral group, which then by some Morita equivalence
leads to a construction of the DaHa for the spin symmetric group.
The twisted group algebra of this double covering group is a
semidirect product $\Cl_n \rtimes \C S_n$ between the group
algebra of the symmetric group $S_n$ and a Clifford algebra
$\Cl_n$ in $n$ generators.
\subsection{}

We introduce in Section~\ref{sec:clifford} the rational {double
affine Hecke-Clifford algebra} (DaHCa) $\An$ with parameter $u$
and show that $\An$ has a standard PBW basis. In particular, $\An$
admits a triangular decomposition with the middle part being the
algebra $\Cl_n \rtimes \C S_n$. In contrast to the usual rational
DaHa, the symmetry between the two halves of polynomial generators
for $\An$ has to be broken. We further provide a realization of
$\An$ in terms of Dunkl operators.

A degenerate affine Hecke-Clifford algebra $\Anhat$ was introduced
by Nazarov \cite{Naz} for studying the Young symmetrizer for the
spin symmetric group (cf. Brundan-Kleshchev \cite{BK, Kle} for
recent development on representations of $\Anhat$ and see
\cite{Dr} and Lusztig \cite{Lu} for the usual degenerate/graded
affine Hecke algebras). We show that the DaHa $\An$ contains a
family of subalgebras isomorphic to $\Anhat$, analogous to the
usual DaHa setup (cf. \cite{Ch, EG}).

We then introduce the rational {spin double affine Hecke algebra}
(sDaHa) $\Bn$ associated to the spin symmetric group and parameter
$u$. The PBW theorem for $\Bn$ provides the following isomorphism
of vector spaces
$$\Bn \cong \Cl[{\xi}_1,\ldots, {\xi}_n]
 \otimes \Tn   \otimes \C[y_1,\ldots, y_n],
$$
where $\Tn$ is the spin symmetric group algebra. A key new feature
is that the generators $y_i$'s commute with each other but the
$\xi_i$'s {\em anti-commute} with each other (cf.
Definition~\ref{def:spindaha}). We also introduce a (new)
degenerate spin affine Hecke algebra $\Bnhat$.
%
The results for $\Bn$ and $\Bnhat$ in Section~\ref{sec:spin} are
parallel to those for $\An$ and $\Anhat$ in
Section~\ref{sec:clifford}.

In Section~\ref{sec:isom} we establish an explicit isomorphism of super algebras
\begin{eqnarray*}
\An \stackrel{\cong}{\longrightarrow} \Cl_n \otimes \Bn.
\end{eqnarray*}
In this case, the superalgebras $\An$ and $\Bn$ are said to be
Morita super-equivalent (cf. \cite[13.2]{Kle} for a justification
of the terminology). We also establish another Morita
super-equivalence between $\Anhat$ and $\Bnhat$ by exhibiting an
superalgebra isomorphism $\Anhat \cong \Cl_n \otimes \Bnhat.$ Such
isomorphisms provide a conceptual explanation and easy proofs for
the results in Section~\ref{sec:spin} parallel to
Section~\ref{sec:clifford}. Without such a connection to $\An$,
even the definition of $\Bn$ is by no means obvious. A
finite-dimensional version of such isomorphisms, discovered
independently in Sergeev \cite{Ser} and Yamaguchi \cite{Yam},
explained for the well-known fact that the representation theory
of the spin symmetric group is (essentially) equivalent to that of
the algebra $\Cl_n \rtimes S_n$. The three isomorphisms in the
finite, affine, and double affine setups are compatible with each
other.

We finally introduce the trigonometric DaHCa $\An_{tr}$ and the
trigonometric sDaHa $\Bn_{tr}$. We construct explicitly a
superalgebra isomorphism
$$\An_{tr} \cong \Cl_n \otimes \Bn_{tr}.
$$
We further establish the precise connections between our rational
and trigonometric DaHa's, following Suzuki \cite{Suz} in the usual
DaHa setup. This is developed in Section~\ref{sec:trig}.
\subsection{}
This paper is a first step of a program in which we attempt to
develop a new theory of spin Hecke algebras.
The new features of the constructions in this paper will be
instrumental in telling us how to relax the usual setup (cf.
\cite{Dr, EG} and Ram-Shepler \cite{RS} for the degenerate case)
in order to define appropriate spin (degenerate/double) affine
Hecke algebras more generally. The algebras $\An$ and $\Bn$ have
large centers, and are expected to have rich finite-dimensional
representation theory (compare \cite{EG} and Gordon \cite{Gor}).
\section{The rational double affine Hecke-Clifford algebras} \label{sec:clifford}
\subsection{The symmetric group $S_n$}

Recall that the symmetric group $S_n$ is generated by $s_i$ $(1 \le
i \le n-1)$ subject to the relations:
\begin{eqnarray}
s_i s_j &=& s_j s_i,\quad |i-j|>1    \nonumber \\
s_i^2 =1,\quad s_i s_{i+1} s_i &=& s_{i+1} s_i s_{i+1}.
\label{symmetric}
\end{eqnarray}
Let $s_{ij}$ denote the transposition of $i$ and $j$ in the
symmetric group $S_n$.
\subsection{The rational DaHCa}

\begin{definition}
Let $u\in \C$. The rational {\em double affine Hecke-Clifford
algebra} (DaHCa) $\An$ is the $\C$-algebra  generated by $x_i,
y_i, c_i (1\le i \le n)$ and $S_n$, subject to the following
relations:
\begin{eqnarray}  \label{poly}
x_i x_j = x_jx_i, && y_i y_j = y_jy_i, \quad (\forall i, j)
  \nonumber \\
\sigma x_i = x_{\sigma i} \sigma, && \sigma y_i = y_{\sigma i}
\sigma  \quad (\sigma \in S_n)  \nonumber  \\
c_i x_i = -x_i c_i, &&c_i y_i =  y_i c_i,  \label{nazarov}  \\
c_j x_i = x_i c_j, && c_j y_i = y_i c_j, \quad (i \neq j)
\nonumber
\end{eqnarray}
\begin{align}
\sigma c_i =c_{\sigma i} \sigma   \; (\sigma \in S_n), \quad
c_i^2 =1, \quad c_i c_j = -c_j c_i  \; (i \neq j)  \label{clifford}
\end{align}
\begin{eqnarray}
[y_j, x_i] &=&
 u(1 +c_j c_i)s_{ij}, \quad (i \neq j) \label{xyij}
\\
{[}y_i, x_i]
  &=&  -u\sum_{k \ne i} (1 +c_kc_i) s_{ki}.   \label{xyii}
\end{eqnarray}
\end{definition}

Denote by $\Cl_n$, or $\Cl (c_1, \ldots, c_n)$, the Clifford
algebra generated by $c_1, \ldots, c_n$. The rational DaHCa $\An$
is a super (i.e. $\Z_2$-graded) algebra with $|c_i| =\bar{1}$ and
$|x_i| = |y_i| = |s_{ij}| =\bar{0}$. Unlike the usual rational
DaHa, the symmetry between $x$'s and $y$'s is broken in the
definition of $\An$, cf.~(\ref{nazarov}).

For $u =0$, we have $\An_{|_{u=0}}= (\C[\underline{x},
\underline{y}] \tilde{\otimes} \Cl_n) \# S_n$, where the tilde
refers to the unusual sign in (\ref{nazarov}). The $\An_{|_{u=0}}$
clearly has the PBW property, i.e. the $\An_{|_{u=0}}$ has a
linear basis $\{\underline{x}^{\underline{a}}
\underline{y}^{\underline{b}} \sigma
\underline{c}^{\underline{\epsilon}}\}$, where $\sigma \in S_n$,
$\underline{a}, \underline{b} \in \Z_+^n$, $\underline{\epsilon}
\in \{0,1 \}^n$ and $\underline{x}^{\underline{a}}$ denotes the
monomials $x_1^{a_1} \cdots x_n^{a_n}$ etc. The algebras $\An$ are
isomorphic for all $u \neq 0$.

\begin{theorem} \label{th:PBW}
The PBW property holds for $\An$. That is, the multiplication of
the subalgebra gives rise to an isomorphism of vector spaces:
 $$\C[x_1,\ldots,x_n]
 \otimes \C S_n \otimes \Cl_n  \otimes \C[y_1,\ldots, y_n]
 \stackrel{\cong}{\longrightarrow} \An. $$
\end{theorem}
\begin{proof}
This can be proved in a similar way as in \cite[Proof of Th.~1.3,
pp.256-7]{EG} (and see {\em loc. cit.} for earlier works on PBW
algebras), with one crucial modification. Set $V =\C^{2n}$ with the
standard coordinates $\{x_i, y_i\}$, and note that $K =\Cl_n \rtimes
\C S_n$ is a semisimple algebra. {\em The new key observation is
that $E :=V\otimes_\C K$ is naturally a $K$-bimodule even though $K$
does not act on $V$.} The right action of $K$ is via the right
multiplication $g: v \otimes a \mapsto v \otimes (ag)$. The $S_n$
acts on the left by $\sigma: v \otimes a \mapsto v^\sigma \otimes
(\sigma a)$, while $\Cl_n$ acts by $c_i: x_j \otimes a \mapsto
(-1)^{\delta_{ij}} x_j \otimes (c_i a)$ and $c_i:  y_j \otimes a
\mapsto y_j \otimes (c_i a)$. The rest of the proof is the same as
in {\em loc. cit.}
\end{proof}

\begin{remark}
The observation in the above proof, when applied to more general
pairs $(V, K)$, leads to a generalization of Drinfeld's setup
\cite{Dr}. The proof above also shows that $\An$ is the only
`quadratic' deformation of $\An_{|_{u=0}}$ with PBW property. Note
the special feature of $\An$ that an additional term $t \cdot 1$
which appeared in  the symplectic reflection algebra \cite{EG} is
not allowed in (\ref{xyii}) by the conjugation invariance with
respect to $c_i$.
\end{remark}

\begin{remark}
Introducing an additional central element $z$ such that $z^2 =1$
to $\An$, we can define a modified algebra $\ddot{\mathfrak
H}^{\text{mod}}$ with relations
$$
c_i c_j = z c_j c_i, \quad c_i x_i = z x_i c_i$$
replacing the corresponding relations in (\ref{clifford}) and
(\ref{nazarov}). Then $\ddot{\mathfrak H}^{\text{mod}} / \langle
z+1 \rangle \cong \An$. On the other hand, $\ddot{\mathfrak
H}^{\text{mod}}/ \langle z-1 \rangle$ is isomorphic to the
ordinary rational DaHa (cf. \cite{EG}) associated to the
reflection group $\Z_2^n \rtimes S_n$ with specialized parameters.
A similar remark applies to the algebra $\Bn$ introduced in
Section~\ref{sec:spin}.
\end{remark}
{Below we always assume that $u \neq 0$ unless otherwise
specified.}
%



%
%
\subsection{The degenerate affine Hecke-Clifford algebra}
\label{sec:heckeclifford}

The degenerate affine Hecke-Clifford algebra was introduced by
Nazarov (also called affine Sergeev algebra in \cite{Naz}). It is
the algebra $\Anhat$ generated by $a_i, c_i (i =1,\ldots, n)$ and
$S_n$, subject to the relations (\ref{symmetric}), (\ref{clifford})
and the following relations:
\begin{eqnarray}
a_j s_i &=& s_i a_j \quad (j \neq i,i+1),  \nonumber \\
a_{i+1} s_i -s_i a_i &=& 1-c_{i+1}c_i,  \nonumber \\
c_i a_i &=& -a_i c_i,  \label{affineNaz} \\
a_i a_j = a_ja_i, &&  c_j a_i = a_i c_j, \quad (i \neq j).
\nonumber
\end{eqnarray}
Our convention $c_i^2 =1$ is consistent with \cite{Kle}, not with
\cite{Naz}. The $\Anhat$ is a superalgebra with $c_i$ being
odd and $a_i$ being even for each $i$.
For $1 \le i \le n$, set
 $$M_i := \sum_{k<i} (1 -c_ic_k) s_{ki}, $$
The $M_i$'s are the Jucys-Murphy elements introduced in \cite{Naz}
and they satisfy
\begin{equation} \label{JMcomm}
M_i M_j =M_j M_i, \quad (\forall i,j).
\end{equation}

\begin{proposition} \cite{Naz}
\begin{enumerate}
 \item
The algebra $\Anhat$ admits the PBW property. That is, the
multiplication of the subalgebras induces a vector space
isomorphism
$$ \C  [a_1, \ldots, a_n] \otimes   \C S_n \otimes  \Cl_n
\stackrel{\cong}{\longrightarrow} \Anhat.$$
 \item
There exists a unique algebra homomorphism $\Anhat \rightarrow
\Cl_n \rtimes \C S_n$, which restricts to the identity map on the
subalgebra $\Cl_n \rtimes \C S_n$ of $\Anhat$, and sends $a_1$ to
$0$. Moreover, this homomorphism sends each $a_i$ to $M_i\;  (1\le
i \le n)$. \item The even center of the algebra $\Anhat$ is $\C
[a_1^2, \ldots, a_n^2]^{S_n}$.
\end{enumerate}
\end{proposition}

The intertwining elements
$$\phi_i :=s_i (a_i^2 -a_{i+1}^2) +(a_i +a_{i+1})
+c_ic_{i+1} (a_i -a_{i+1}), \;\; 1\le i \le n-1$$
were introduced in \cite{Naz} (cf. \cite[Chapter~14]{Kle}).

\begin{proposition}  \cite{Naz}
The intertwining elements $\phi_i$'s satisfy
\begin{eqnarray*}
\phi_i^2 &=& 2a_i^2 +2a_{i+1}^2 - (a_i^2 -a_{i+1}^2)^2 \\
\phi_i \phi_{i+1} \phi_i &=& \phi_{i+1} \phi_i \phi_{i+1} \\
\phi_i\phi_j &=& \phi_j \phi_i \qquad (|i-j|>1).
\end{eqnarray*}
\end{proposition}
\subsection{The affine Hecke-Clifford subalgebra of $\An$} \label{sec:affinesubalg}

Set
$$z_i := u^{-1} y_i x_i + M_i.$$

\begin{lemma} \label{commfamily}
For all $i,j$ and each $\alpha \in \C$, the following identities
hold:
\begin{eqnarray}
 [z_i, z_j] &=&0, \label{eq:comm1}  \\
 {[}x_i, z_j] -[x_j, z_i] &=&0, \label{eq:comm2}  \\
 {[}\alpha x_i +z_i, \alpha x_j +z_j] &=&0. \label{eq:comm3}
\end{eqnarray}
\end{lemma}

\begin{proof}
We may assume $i<j$. Then, by (\ref{JMcomm}), we have
\begin{eqnarray*}
 u[z_i, z_j]
 &=& u^{-1} [y_ix_i, y_jx_j] +[y_ix_i, M_j] \\
 &=& u^{-1}  (y_i [x_i, y_j] x_j +y_j [y_i, x_j]x_i)
 + [y_ix_i, s_{ij} -c_j c_is_{ij}]
 \\
 &=& -y_i (s_{ij} -c_ic_js_{ij}) x_j
 +y_j (s_{ij} +c_ic_j s_{ij})x_i \\
 &&
  + (y_ix_i -y_jx_j) s_{ij} - (y_ix_i +y_jx_j) c_j c_i s_{ij}
 =0.
\end{eqnarray*}

To prove (\ref{eq:comm2}) for $i<j$, we calculate that
\begin{eqnarray*}
 && [x_i, z_j] -[x_j, z_i]\\
 &=& [x_i, u^{-1}  y_jx_j + M_j] -   [x_j, u^{-1} y_ix_i] \\
 &=& u^{-1} [x_i, y_j]x_j +[x_i, s_{ij}-c_jc_is_{ij}] -u^{-1}
 [x_j, y_i] x_i \\
 &=& (c_ic_js_{ij}-s_{ij}) x_j +(x_i -x_j)s_{ij}
 -(x_i +x_j) c_jc_i s_{ij} +  (s_{ij} +c_ic_j s_{ij}) x_i
 = 0.
\end{eqnarray*}
Now (\ref{eq:comm3}) follows from (\ref{eq:comm1}) and
(\ref{eq:comm2}).
\end{proof}

We remark that $[y_i, z_j] - [y_j, z_i] \neq 0 \quad (i \neq j)$.

\begin{lemma} \label{lem:morecomm}
The following identities hold:
$$c_i z_i = -z_i c_i, \quad c_i z_j =z_j c_i \quad (i \neq j).$$
\end{lemma}

\begin{proof}
Recall $z_i = u^{-1} y_ix_i + M_i$. It is known \cite{Naz} that
the Jucys-Murphy elements satisfy  $ c_i M_i = -M_i c_i$, and $c_j
M_i = M_i c_j$ for $ i \neq j.$ Clearly, $c_i (y_ix_i) = -(y_ix_i)
c_i,$ and $c_j (y_ix_i) = (y_ix_i) c_j$ for $ i \neq j.$ Now the
lemma follows.
\end{proof}

\begin{lemma} \label{lem:heckecomm}
The following identities hold:
\begin{eqnarray*}
 z_{i+1} s_i - s_i z_i &=& 1-c_{i+1} c_i, \\
(\alpha x_{i+1} +z_{i+1}) s_i -s_i (\alpha x_i +z_i)  &=&
1-c_{i+1} c_i.
\end{eqnarray*}
\end{lemma}
\begin{proof}
The second identity follows from the first one and the equation
$x_{i+1} s_i =s_i x_i$. Note that
\begin{eqnarray*}
z_{i+1} s_i - s_i z_i
 &=& (u^{-1} y_{i+1} x_{i+1} +M_{i+1}) s_i -s_i (u^{-1} y_ix_i +M_i)\\
 &=& M_{i+1}s_i -s_i M_i
 = 1-c_{i+1} c_i.
\end{eqnarray*}
This proves the first identity.
\end{proof}

The next theorem follows now from Lemmas~\ref{commfamily},
\ref{lem:morecomm},  and \ref{lem:heckecomm}.
\begin{theorem} \label{th:affinealgA}
Fix $\alpha \in \C$. The subalgebra of $\An$ generated by $c_i,
\alpha x_i + z_i$ $(1\le i \le n)$ and $S_n$ is isomorphic to the
degenerate affine Hecke-Clifford algebra $\Anhat$.
\end{theorem}

\subsection{The Dunkl operator for $\An$}
The PBW Theorem~\ref{th:PBW} provides the triangular decomposition
of $\An$:
$$\An \cong \C[\underline{x}] \otimes
(\Cl_n \rtimes \C S_n) \otimes \C[\underline{y}]$$
where $\C[\underline{x}]$ is the short-hand for $\C[x_1, \ldots,
x_n]$, etc. Denote by $\Anhat_x$ the subalgebra of $\An$ generated
by $\Cl_n \rtimes \C S_n$ and $x_1, \ldots, x_n$. The subalgebra
$\Anhat_y$ is similarly defined. Take any $\Cl_n \rtimes \C
S_n$-module $W$, and extend it to an $\Anhat_x$-module with the
trivial action of $x_i$'s. Consider the induced $\An$-module
$\ind^{\An}_{\Anhat_x} {W}$, which as a vector space is isomorphic
to $\C[y_1, \ldots, y_n] \otimes W.$ Denote the action by $\circ$.
The most interesting $S_n \ltimes \Cl_n$-module is the so-called
{\em basic spin module} $L_n =\Cl (c_1, \ldots, c_n)$ of $S_n
\ltimes \Cl_n$, (where $\Cl_n$ acts by left multiplication and
$S_n$ acts by permuting the $c_i$'s), and the induced $\An$-module
$\ind^{\An}_{\Anhat_x} {L_n}$ is then identified as $\C[y_1,
\ldots, y_n] \otimes \Cl (c_1, \ldots, c_n).$
Given operators $h$ and $g$, $\frac{h}g$ stands for $g^{-1}h$
throughout the paper.
\begin{theorem}  \label{th:dunkl}
Given a $\Cl_n \rtimes \C S_n$-module $W$, the action of $x_i$ on
$\C[y_1, \ldots, y_n] \otimes W$ is realized as ``Dunkl operator" as
follows. For any polynomial $f =f(\underline{y})$ and $w \in W$, we
have
\begin{eqnarray*}
 x_i \circ (f \otimes w ) =
  u \sum_{k \neq i} \frac{(1 -s_{ki})(f)}{y_i -y_k} \otimes
  (1 -c_ic_k) s_{ki} (w).
\end{eqnarray*}
\end{theorem}

\begin{proof}
We calculate that
\begin{eqnarray*}
 x_i \circ (f \otimes w )
 &=&  [x_i, f] \circ w + f x_i \circ w
 =  [x_i, f] \circ w.
\end{eqnarray*}
Now the theorem follows from Lemma~\ref{lem:dunkl}~(2) below.
\end{proof}

\begin{lemma} \label{lem:dunkl}
\begin{enumerate}
 \item
For $a \in \N$ and $i \neq j$, we have
\begin{eqnarray*}
[x_i, y_j^a]
 &=& - u  \frac{y_i^a -y_j^a}{y_i -y_j} (1 -c_ic_j) s_{ij}
\\
{[} x_i, y_i^a]
 &=& u \sum_{k \neq i}
  \frac{y_i^a -y_k^a}{y_i -y_k} (1 -c_ic_k) s_{ki}.
\end{eqnarray*}
 \item
Let $f(\underline{y})$ be any polynomial in $y_1, \ldots, y_n$.
Then,
\begin{eqnarray*}
 [x_i, f(\underline{y})]
 = u \sum_{k \neq i}
  \frac{(1 -s_{ki})f(\underline{y})}{y_i -y_k} (1 -c_ic_k) s_{ki}.
\end{eqnarray*}
\end{enumerate}
\end{lemma}
\begin{proof}
(1) is proved by induction on $a$ using $[x_i, y_j^a]
=[x_i,y_j^{a-1}]y_j +y_j^{a-1} [x_i, y_j]$.

(2) It suffices to check for every monomial $f$. We first prove
the case when $f =\prod_{j \neq i} y_j^{a_j}$ using the first
identity in (1). Then we combine with the second identity in (1)
to prove for a general monomial $f$ involving powers of $y_i$.
\end{proof}

Note that the symmetry between $x$'s and $y$'s in our algebra
$\An$ is broken. We compute below the Dunkl operator for $y_i$.
\begin{lemma}  \label{lem:commMore}
 For $a \in \N$ and $1 \le i \neq j \le n,$ we have
\begin{eqnarray*}
 [y_i, x_i^a]
 &=& u \sum_{k \neq i} \left (
  -\frac{x_k^a -x_i^a}{x_k -x_i} s_{ki}
  + \frac{x_i^a - (-x_k)^a}{x_k +x_i} c_i s_{ki}c_i
  \right ) \\
 {[}y_i, x_j^a]
 &=& u \frac{x_j^a -x_i^a}{x_j -x_i} s_{ij}
  + u \frac{x_j^a - (-x_i)^a}{x_j +x_i} c_i s_{ij}c_i.
\end{eqnarray*}
\end{lemma}

\begin{proof}
Follows by induction on $a$.
\end{proof}

We can identify $\C[x_1, \cdots, x_n] \otimes W$ with the induced
$\An$-module $\ind^{\An}_{\Anhat_y} W$, where $W$ is a $\C S_n
\ltimes \Cl_n$-module and is then extended to an $\Anhat_y$-module
trivially.

\begin{proposition}
The action of $y_i$ on $\C[x_1, \cdots, x_n] \otimes W$ is given as
follows. For any polynomial $f(\underline{x})$ and $w\in W$, we have
\begin{eqnarray*}
 y_i \circ (f\otimes w) =
  \sum_{k \neq i}
 \left( \frac{f -s_{ki}(f)}{x_k -x_i} \otimes s_{ki} (w)
  + \frac{f c_ic_k - c_ic_k s_{ki} (f)}{x_k +x_i} \otimes s_{ik}(w)
 \right).
\end{eqnarray*}
\end{proposition}

\begin{proof}
First note that $y_i$ commutes with all $c_j$'s. So we are reduced
to the case $f=1$ and $w$ being a monomial. This follows by
induction from Lemma~\ref{lem:commMore}.
%
\end{proof}

\subsection{The center of $\An$}
Denote by $\Zn (A)$ the even center of a superalgebra $A$, which
consists of all the even central elements in $A$.
\begin{proposition} \label{smallcenter}
\begin{enumerate}
 \item
For $u=0$, the even center $\Zn (\An_{|_{u=0}})$ of the
superalgebra $\An_{|_{u=0}}$ is $\C[x_1^2, \ldots, x_n^2; y_1,
\ldots, y_n]^{\Delta S_n}$, where $\Delta$ denotes the diagonal
action on $x_j^2$'s and $y_j$'s.
 \item
We have $\C[x_1^2, \ldots, x_n^2]^{S_n} \subset \Zn (\An)$ and
$\C[y_1, \ldots, y_n]^{S_n} \subset \Zn (\An)$, $\forall u \in
\C$.
\end{enumerate}
\end{proposition}
\begin{proof}
(1) is clear.

(2) By Lemma~\ref{lem:dunkl}~(1), for each $i$ we have that
$$[x_i, y_1^k +\cdots + y_n^k] =0, \quad \forall k \geq 1.$$
By definition of $\An$, $y_1^k +\cdots +y_n^k$ commutes with $\C
S_n \ltimes \Cl_n$. So $y_1^k +\cdots +y_n^k \in  \Zn (\An)$ for
all $k$, and thus $\C[y_1, \ldots, y_n]^{S_n} \subset \Zn (\An)$.
We have $\C[x_1^2, \ldots, x_n^2]^{S_n} \subset \Zn (\An)$, by
arguing similarly using Lemma~\ref{lem:commMore} and the
definition of $\An$.
\end{proof}


\begin{example}
For $n=2$, $x_1^2 y_1 +x_2^2 y_2 -(x_1 +x_2)s_{12}
-c_1(x_1+x_2)s_{12}c_1 \in \Zn (\An)$.
\end{example}


\section{The rational spin double affine Hecke algebras} \label{sec:spin}

The results in this section are presented in a way parallel to
Section~\ref{sec:clifford}, but with all proofs omitted except for the
PBW theorem. They can be
either proved in the same way as for their counterparts in
Section~\ref{sec:clifford}, or follow directly from the
counterparts in Section~\ref{sec:clifford} via the isomorphisms to
be established in Section~\ref{sec:isom}.
\subsection{The spin symmetric group algebra}

The symmetric group $S_n$ affords a double cover $\tS n$,
nontrivial for $n\ge 4$, according to Schur \cite{Sch} (cf.
\cite{Kle}):
\begin{equation*}
1 \longrightarrow\mathbb Z_2  {\longrightarrow} \tS n
 {\longrightarrow} S_n \longrightarrow 1.
\end{equation*}
Denote $\mathbb Z_2 =\{1,z\}$. The {\em spin symmetric group
algebra} $\Tn :=\C[\tSn]/\langle z+1 \rangle$ is generated by
$t_i, i=1, \cdots, n-1$ subject to the following relations:
\begin{eqnarray}\label{E:defrel1}
 t_i^2=1, && t_it_{i+1}t_i =t_{i+1}t_it_{i+1}, \nonumber
\\
 t_it_j=-t_jt_i, &&   |i -j|>1.
\end{eqnarray}
The algebra $\Tn$ is naturally a superalgebra by declaring $t_i$
for every $i$ to be odd.

Define the ``transpositions" of odd parity, for $1 \le i<j \le n$,
\begin{eqnarray*}
[i,j] =- [j,i] = (-1)^{j-i-1} t_{j-1} \cdots t_{i+1} t_i t_{i+1}
\cdots t_{j-1}
\end{eqnarray*}
which satisfies the following relations:
\begin{eqnarray*}
[i, i+1] =t_i,&& [i,j]^2 = 1, \\
t_i[i,j]t_i = -[i+1,j], && \text{for } j \neq i,i+1.
\end{eqnarray*}
Define the {\em odd Jucys-Murphy elements}
\begin{eqnarray*}
\texttt{M}_j =\sum_{k<j} [k,j], \qquad (1 \le j \le n).
\end{eqnarray*}
We have \cite{Ser} (also see \cite[Chapter
13]{Kle})
$$[\texttt{M}_i, \texttt{M}_j]_+ :=
 \texttt{M}_i \texttt{M}_j +\texttt{M}_j \texttt{M}_i =0, \qquad (i
\neq j).$$
As usual, we have denoted that $[a, b]_+ =ab +ba.$

\subsection{The rational sDaHa}

\begin{definition}  \label{def:spindaha}
Let $u\in \C$. The rational {\em spin double affine Hecke algebra}
(sDaHa) is the algebra $\Bn$ generated by ${\xi}_i, y_i (1\le i
\le n)$ and $t_i (1\le i \le n-1)$, subject to the relations
(\ref{E:defrel1}) for $t_i$ and the following relations:
\begin{eqnarray*}
{\xi}_i {\xi}_j = - {\xi}_j {\xi}_i, && y_i
y_j = y_jy_i,  \quad (i \neq j), \label{polyodd1} \\
 t_i {\xi}_i = -{\xi}_{i+1} t_i, && t_i y_i = y_{i+1}
t_i,   \\
t_j {\xi}_i = -{\xi}_{i} t_j, && t_j y_i = y_{i} t_j, \quad (i
\neq j, j+1),
\end{eqnarray*}
\begin{eqnarray*}
[y_i, {\xi}_j] &=&
 u [i,j], \quad (i \neq j), \\
 {[}y_i, {\xi}_i] &=&
   u\sum_{k \ne i} [i,k].
\end{eqnarray*}
%
\end{definition}
Denote by $\Cl[{\xi}_1,\ldots, {\xi}_n]$ the algebra generated by
${\xi}_1,\ldots, {\xi}_n$ subject to the relations ${\xi}_i
{\xi}_j = - {\xi}_j {\xi}_i$ for $i \neq j$. Clearly, it has a
linear basis $\underline{\xi}^{\underline{\texttt{a}}}$, where
$\underline{\texttt{a}} \in \Z_+^n$.
\begin{theorem}   \label{th:PBW-Bn}
 The PBW property holds for $\Bn$. That is,
 the multiplication of
 the subalgebras gives rise to an isomorphism of vector spaces:
 $$\Cl[{\xi}_1,\ldots, {\xi}_n]
 \otimes \Tn   \otimes \C[y_1,\ldots, y_n] \stackrel{\cong}{\longrightarrow} \Bn. $$
\end{theorem}

\begin{proof}
It follows readily from the definition that $\Bn$ is spanned by
the elements $\underline{{\xi}}^{\underline{a}} \sigma
\underline{y}^{\underline{b}}$, where $\sigma$ runs over a basis
for $\Tn$ and $\underline{a}, \underline{b} \in \Z_+^n$. Recall
the homomorphism $\Psi$ in Theorem~\ref{th:isoDAHA}. The images
$\Psi( \underline{{\xi}}^{\underline{a}} \sigma
\underline{y}^{\underline{b}})$ in explicit formulas are clearly
all linearly independent in $\An$ by the PBW Theorem~\ref{th:PBW}
for $\An$. Thus, all the elements
$\underline{{\xi}}^{\underline{a}} \sigma
\underline{y}^{\underline{b}}$ are linearly independent in $\Bn$.
\end{proof}

\begin{remark}
Note the subtle signs in the definition of $\Bn$. The algebra
$\Bn$ is naturally a superalgebra with $y_i$ being even and
${\xi}_i, t_i$ being odd for each $i$. The algebras $\Bn$ are all
isomorphic for different $u \neq 0$.

{\em Below we assume that $u \neq 0$ unless otherwise specified.}
\end{remark}

\subsection{The degenerate spin affine Hecke algebra}
\label{sec:affineHeckeSpin}

The {\em degenerate spin affine Hecke algebra} is the algebra
$\Bnhat$ generated by ${b}_i (1\le i \le n)$ and $t_i (1\le i \le
n-1)$, subject to the relations (\ref{E:defrel1}) for $t_i$ and the
following relations:
\begin{eqnarray*}
{b}_i {b}_j &=& - {b}_j {b}_i,
\quad (i \neq j)  \\
  {b}_{i+1} t_i &=&   -  t_i {b}_i +1,   \\
t_j {b}_i &=& - {b}_{i} t_j,  \quad (i \neq j, j+1).
\end{eqnarray*}
The algebra $\Bnhat$ has a superalgebra structure with both $b_i$
and $t_i$ being odd elements for every $i$. Note that the algebra
$\Bnhat$ contains $\Cl [b_1, \ldots, b_n]$ and $\Tn$ as
subalgebras.

\begin{proposition}
 \begin{enumerate}
 \item
The algebra $\Bnhat$ admits the PBW property. That is, the
multiplication of the subalgebras induces a vector space isomorphism
$$ \Cl [b_1, \ldots, b_n] \otimes \Tn
\stackrel{\cong}{\longrightarrow} \Bnhat.$$
 \item
There exists a unique algebra homomorphism $\Bnhat \rightarrow
\Tn$, which extends the identity map on the subalgebra $\Tn$ of
$\Bnhat$ and sends $b_1$ to $0$. Moreover, this homomorphism sends
each $b_i$ to $\texttt{M}_i \; (1 \le i\le n)$. \item The even
center of the superalgebra $\Bnhat$ is $\C  [b_1^2, \ldots,
b_n^2]^{S_n}$.
\end{enumerate}
\end{proposition}

We define the (odd) intertwining elements in $\Bnhat$:
$$\psi_i := t_i (b_i^2 -b_{i+1}^2) - (b_i -b_{i+1}),\;\; 1\le i \le n-1.$$

\begin{proposition}
The intertwining elements $\psi_i$'s satisfy that
\begin{eqnarray*}
\psi_i^2 &=& b_i^2 +b_{i+1}^2 - (b_i^2 -b_{i+1}^2)^2 \\
\psi_i \psi_{i+1} \psi_i &=& \psi_{i+1} \psi_i \psi_{i+1} \\
\psi_i\psi_j &=& - \psi_j \psi_i \;\; (|i-j|>1) \\
\psi_i b_j &=& -b_j \psi_i \;\;\, (j \neq i,i+1)   \\
\psi_i b_i = -b_{i+1} \psi_i, &&
\psi_i b_{i+1} = -b_{i} \psi_i.
\end{eqnarray*}
\end{proposition}

\subsection{The spin affine Hecke subalgebra of $\Bn$}

Let $$\mathfrak{z}_i := u^{-1}y_i{\xi}_i + \texttt{M}_i.$$

\begin{lemma} \label{supercommfamily}
For $i\neq j$, the following identities hold:
\begin{eqnarray*}
 [\mathfrak{z}_i, \mathfrak{z}_j]_+ &=& 0, \\
{[}{\xi}_i, \mathfrak{z}_j]_+
 +[{\xi}_j, \mathfrak{z}_i]_+ &=& 0,  \\
 {[}\alpha {\xi}_i + \mathfrak{z}_i, \alpha {\xi}_j +
 \mathfrak{z}_j]_+  &=& 0, \quad \forall \alpha \in \C.
\end{eqnarray*}
\end{lemma}

\begin{lemma}
The following identities hold:
\begin{eqnarray*}
 %
(\alpha {\xi}_{i+1} + \mathfrak{z}_{i+1}) t_i
 &=& - t_i (\alpha {\xi}_i + \mathfrak{z}_i) +1,
 \quad \forall \alpha \in \C.
\end{eqnarray*}
\end{lemma}

\begin{theorem} \label{th:affinealgB}
Fix $\alpha \in \C$. The subalgebra of $\Bn$ generated by $ \alpha
{\xi}_i + \mathfrak{z}_i$ $(1\le i \le n)$ and $\Tn$ is isomorphic
to the degenerate spin affine Hecke algebra $\Bnhat$.
\end{theorem}

\subsection{The Dunkl operators for $\Bn$}

Denote by $\Bnhat_\xi$ the subalgebra of $\Bn$ generated by $\Tn$
and ${\xi}_1, \ldots, {\xi}_n$. Take any $\Tn$-module $W$, and
extend it to an $\Bnhat_\xi$-module with the trivial action of
${\xi}_i$'s. Consider the induced $\Bn$-module
$\ind^{\Bn}_{\Bnhat_\xi} W  \cong \C[y_1, \ldots, y_n] \otimes W.$

\begin{theorem}
The action of $\xi_i \in \Bn$ on $\C[y_1, \ldots, y_n] \otimes W$
is realized as ``Dunkl operator" as follows. For any polynomial $f
=f(\underline{y})$ and $w \in W$, we have
\begin{eqnarray*}
 {\xi}_i \circ (f \otimes w) =
 u \sum_{k \neq i} \frac{(1 -s_{ki})(f)}{y_i -y_k}
  \otimes [k,i](w).
\end{eqnarray*}
\end{theorem}

\subsection{The center of $\Bn$}

\begin{proposition}
\begin{enumerate}
 \item
For $u=0$, the even center $\mathcal Z(\Bn_{|_{u=0}})$ of the
superalgebra $\Bn_{|_{u=0}}$ is $ \C[{\xi}_1^2, \ldots, {\xi}_n^2;
y_1, \ldots, y_n]^{\Delta S_n}$, where $\Delta$ denotes the
diagonal action on the $\xi_j^2$'s and $y_j$'s.
 \item
 For any $u\in \C$, $\C[{\xi}_1^2, \ldots, {\xi}_n^2]^{S_n}
\subset \Zn (\Bn)$ and $\C[y_1, \ldots, y_n]^{S_n} \subset \Zn
(\Bn)$.
\end{enumerate}
\end{proposition}

\begin{example}
 For $n=2$, ${\xi}_1^2 y_1 +{\xi}_2^2 y_2
 +2 ({\xi}_1-{\xi}_2)t_1 \in \Zn (\Bn)$.
\end{example}

\section{The algebra isomorphisms} \label{sec:isom}
The results in this section use earlier notations but their proofs
will be independent of Section~\ref{sec:spin}. Note that the
multiplication in a tensor product $\Cl \otimes \Bn$ of two super
algebras $\Cl$ and $\Bn$ has a suitable sign convention:
$$(c' \otimes b')(c \otimes b) =(-1)^{|b'||c|} (c'c
\otimes b'b).$$
 We shall write a typical element in $\Cl \otimes
\Bn$ as $cb$ rather than $c \otimes b$, and use short-hand notations
$c =c \otimes 1, b =1\otimes b$.
\begin{theorem}
There exists a superalgebra isomorphism:
\begin{align*}
 \widehat{\Phi}: &\;
  \Anhat \stackrel{\cong}{\longrightarrow} \Cl_n \bigotimes
  \Bnhat, \\
c_i \mapsto c_i,   \; a_i \mapsto &  \sqrt{-2} c_i  b_i,  \;
 s_i \mapsto \frac{1}{\sqrt{-2}} (c_i -c_{i+1})  t_i.
\end{align*}
The inverse map is given by
\begin{eqnarray*}
 \widehat{\Psi}: c_i \mapsto c_i,  \;
 b_i \mapsto \frac{1}{\sqrt{-2}}c_i a_i, \;
 t_i \mapsto  \frac{1}{\sqrt{-2}} (c_{i+1} -c_i) s_i.
\end{eqnarray*}
\end{theorem}
\begin{proof}
Clearly $\widehat{\Phi},  \widehat{\Psi}$ are inverses to each
other, once we show that they are algebra homomorphisms. Note that
$\Cl_n\rtimes \C S_n$ is a subalgebra of $\Anhat$ and
$\widehat{\Phi}$ extends the isomorphism established in \cite{Ser,
Yam}
$$\Phi^{fin}: \Cl_n\rtimes \C S_n \longrightarrow
 \Cl_n \otimes \Tn. $$
This takes care of all the defining relations for the images of
$\widehat{\Phi}$ or $\widehat{\Psi}$ of the generators which do
not involve $a_i$ or $b_i$ (the verification of these relations
can be done directly without difficulty). Among the defining
relations involving $b_i$, the most complicated one is
$\widehat{\Psi} ({b}_{i+1} t_i +t_i {b}_i) = 1$, which is
equivalent to the following identity in $\Anhat$:
$$ \frac{1}{\sqrt{-2}}c_{i+1} a_{i+1}
 \cdot \frac{1}{\sqrt{-2}} (c_{i+1} -c_i) s_i
 + \frac{1}{\sqrt{-2}} (c_{i+1} -c_i) s_i
 \cdot \frac{1}{\sqrt{-2}}c_ia_i =1.$$
 This can be checked directly using the defining relations (\ref{affineNaz}) of $\Anhat$.
%
The verification of the remaining relations is similar, and will
be skipped.
\end{proof}
An argument similar to the above implies the following isomorphism
theorem.
\begin{theorem} \label{th:isoDAHA}
There exists a superalgebra isomorphism:
\begin{align*}
 \Phi: &\; \An \stackrel{\cong}{\longrightarrow} \Cl_n \bigotimes
 \Bn,
 \\
y_i \mapsto  y_i, \; c_i \mapsto c_i, & \;
 x_i \mapsto {\sqrt{-2}}c_i  {\xi}_i, \;
 s_i \mapsto \frac{1}{\sqrt{-2}} (c_i -c_{i+1})  t_i.
\end{align*}
The inverse map is given by
\begin{eqnarray*}
\Psi: y_i \mapsto y_i, \; c_i  \mapsto c_i, \;
  {\xi}_i \mapsto \frac{1}{\sqrt{-2}} c_i x_i, \;
 t_i \mapsto \frac{1}{\sqrt{-2}} (c_{i+1}-c_i) s_i.
\end{eqnarray*}
\end{theorem}

\begin{remark}
\begin{enumerate}
\item Recall the intertwining elements $\phi_i$ and $\psi_i$ from
Sections~\ref{sec:clifford} and \ref{sec:spin}. The isomorphism
$\widehat{\Phi}$ sends $\phi_i$ to $-\sqrt{-2} (c_i -c_{i+1})
\psi_i$ for every $i$. \item Recall from
Propositions~\ref{th:affinealgA} and \ref{th:affinealgB} that
$\alpha {x}_i +  {z}_i$ (resp. $\alpha {\xi}_i + \mathfrak{z}_i$)
are the `polynomial' generators of a degenerate affine Hecke
algebra isomorphic to $\Anhat$ (resp. $\Bnhat$). The $\Phi$ sends
$\alpha {x}_i + {z}_i$ to $\sqrt{-2} c_i(\alpha {\xi}_i +
\mathfrak{z}_i)$ for each $i$. This is compatible with the fact
that $\widehat{\Phi}$ sends $a_i$ to $\sqrt{-2} c_ib_i$. \item
Note also that $\widehat{\Phi}$ or $\Phi$ sends the JM element
$M_i$ to $\sqrt{-2} c_i \texttt{M}_i$, which amounts to checking
the identity
$\Phi (\frac1{ \sqrt{-2}}(c_i -c_k)s_{ik}) = [k,i].$
\end{enumerate}
\end{remark}

\section{The trigonometric version of DaHCa and sDaHa} \label{sec:trig}
\subsection{The trigonometric DaHCa}

\begin{definition} \label{def:trigDAHA}
 {\em The trigonometric double affine Hecke-Clifford algebra} is
 the algebra $\An_{tr}$ generated by $e^{\pm \epsilon_i},
 \epsilon_i^\vee, c_i$ $(1 \le i \le n)$ and $S_n$, subject to the following
 relations:
 \begin{eqnarray*}
 c_i^2 =1, && c_i c_j = -c_j c_i  \quad (i \neq j) \\
 e^{\epsilon_i} e^{-\epsilon_i}  =1,
 &&
 e^{\epsilon_i} e^{\epsilon_j} = e^{\epsilon_j} e^{\epsilon_i}, \\
  \epsilon_i^\vee  \epsilon_j^\vee =  \epsilon_j^\vee  \epsilon_i^\vee,
  && c_j e^{\epsilon_i} =e^{\epsilon_i} c_j  \quad (\forall i, j), \\
 \sigma e^{\epsilon_i} =e^{\epsilon_{\sigma i}} \sigma,
 &&
 \sigma c_i =c_{\sigma i} \sigma   \quad (\sigma \in S_n), \\
 c_i \epsilon_i^\vee =-\epsilon_i^\vee c_i,
 &&
 c_j \epsilon_i^\vee = \epsilon_i^\vee c_j \quad (i \neq j), \\
 \epsilon^\vee_{i+1} s_i -s_i \epsilon_i^\vee = u(1-c_{i+1}c_i),
 &&
 \epsilon^\vee_j s_i =s_i \epsilon_j^\vee \quad (j \neq i, i+1), \\
 {[}\epsilon_i^\vee, e^{\eta}]
  &=& u \sum_{k \neq i}  {\text{sgn}(k-i)}
 \frac{e^\eta -e^{s_{ki} (\eta)}}{1 -e^{\text{sgn}(k-i) \cdot (\epsilon_k -\epsilon_i)}}
 (1 -c_i c_k) s_{ki}.
 \end{eqnarray*}
\end{definition}
The algebra $\An_{tr}$ admits a natural superalgebra structure
with $c_i$ being odd and all other generators being even.
 Note that the subalgebra generated by $e^{\epsilon_i}$
$(1 \le i \le n)$, denoted by $\C [P]$, is identified with the group
algebra of the weight lattice of type $GL_n$; the subalgebra
generated by $e^{\epsilon_i}$ $(1 \le i \le n)$ and $S_n$ is
identified with the group algebra of the extended affine Weyl group
of type $GL_n$; for $u\neq 0$, the subalgebra generated by $u^{-1}
\epsilon_i^\vee, c_i$ $(1 \le i \le n)$ and $S_n$ is identified with
the degenerate affine Hecke-Clifford algebra (cf.
Sect.~\ref{sec:heckeclifford}).

\begin{theorem}  \label{th:rat=tr}
\begin{enumerate}
 \item
There exists a unique superalgebras homomorphism $\iota: \An
\rightarrow \An_{tr}$, which extends the identity map on the
subalgebra $\Cl_n \rtimes \C S_n$ such that
 \begin{eqnarray*}
 \iota (y_i) &=& e^{\epsilon_i} \\
 \iota (x_i) &=& e^{-\epsilon_i} \left(\epsilon_i^\vee -u \sum_{k<i}
 (1 -c_i c_k) s_{ki} \right).
 \end{eqnarray*}
 Moreover, $\iota$ is injective.
 \item
 The homomorphism $\iota$ extends to an algebra isomorphism
$$\C[\underline{y}^\pm] \otimes_{\C [ \underline{y} ]} \An
\stackrel{\cong}{\longrightarrow} \An_{tr}.$$ The inverse $j$ is
an extension of the identity map on $\Cl_n \rtimes \C S_n$  such
that
 \begin{eqnarray*}
 j (e^{\epsilon_i}) &=& y_i,  \\
 j (\epsilon_i^\vee) &=& y_i x_i + u\sum_{k<i} (1 -c_ic_k) s_{ki}.
 \end{eqnarray*}
\end{enumerate}
\end{theorem}

\begin{proof}
The $j$ and (the extended) $\iota$ are clearly inverses to each
other, once we show that they are algebra homomorphisms.

Recall from (\ref{sec:affinesubalg}) that $uz_i =y_i x_i + u
\sum_{k<i} (1 -c_ic_k) s_{ki}.$ Let us verify
$$j
([\epsilon_i^\vee, e^{\epsilon_i}]) =[j (\epsilon_i^\vee), j
(e^{\epsilon_i})] =[uz_i, y_i].$$ By a long but straightforward
computation we have that
\begin{eqnarray*}
 {[}uz_i, y_j] &=& -uy_i (1 -c_ic_j) s_{ij}, \quad (i<j), \\
 {[}uz_i, y_j] &=& -uy_j (1-c_ic_j) s_{ij}, \quad (i>j), \\
 {[}uz_i, y_i] &=& u\sum_{k>i} y_i (1 -c_ic_k) s_{ki} + u \sum_{k<i} (1
 -c_ic_k) s_{ki}y_i.
\end{eqnarray*}
By another direct computation, we see that the commutator
${[}\epsilon_i^\vee, e^{\epsilon_j}]$ are given by exactly the
right-hand side of ${[}uz_i, y_j]$ above with $y_i$ therein
replaced by $e^{\epsilon_i}$, for $i<j,$ $i=j$, or $i>j$.

The verification of other relations are even easier and will be
skipped.
\end{proof}

\begin{remark}
Theorem~\ref{th:rat=tr} and Theorem~\ref{th:rat=trBn} below are
analogous to the relation between the usual rational DaHa and
trigonometric DaHa due to Suzuki \cite{Suz}.
\end{remark}


\begin{proposition}
The PBW property holds for $\An_{tr}$. That is, the multiplication
of the subalgebras gives rise to an isomorphism of vector spaces:
$$\C[P] \otimes \C S_n \otimes \Cl_n \otimes \C[ \epsilon_1^\vee,
\ldots,  \epsilon_n^\vee] \stackrel{\cong}{\longrightarrow}
\An_{tr}.
$$
\end{proposition}

\begin{proof}
Clearly the multiplication homomorphism is surjective.

The injectivity of this homomorphism follows from
Theorem~\ref{th:rat=tr}~(2) and the PBW property for $\An$ (and
thus for $\C[\underline{y}^\pm] \otimes_{\C [ \underline{y} ]}
\An$). More explicitly, $j$ sends ${\epsilon_i}^{\vee}$ to $y_i
x_i$ plus lower terms. Via $j$ and with such lower terms ignored,
any linear relation among
${\underline{\epsilon}^{\vee}}^{\underline{a}} \sigma
\underline{c}^{\underline{e}} \underline{y}^{\underline{b}}$ is
translated to the same relation among
$(\underline{yx})^{\underline{a}} \sigma
\underline{c}^{\underline{e}} \underline{y}^{\underline{b}}$. But
we know the latter are all linearly independent.
\end{proof}
\subsection{The trigonometric sDaHa}

\begin{definition}
 The {\em trigonometric spin double affine Hecke algebra} is
 the algebra $\Bn_{tr}$ generated by $e^{\pm \epsilon_i},
 {\zeta}_i$ $(1 \le i \le n)$ and $t_i (1 \le i \le n-1)$, subject to
 the relations (\ref{E:defrel1}) for $t_i$ and the following relations:
 \begin{eqnarray*}
 e^{\epsilon_i} e^{-\epsilon_i}  =1,
 &&
 e^{\epsilon_i} e^{\epsilon_j} = e^{\epsilon_j} e^{\epsilon_i}, \\
 t_i e^{\epsilon_i} =e^{\epsilon_{i+1}} t_i, &&
  t_i e^{\epsilon_j} =e^{\epsilon_{j}} t_i\;\;\; (j \neq i,i+1), \\
  {\zeta}_{i+1} t_i +t_i {\zeta}_i = u,
  &&
  {\zeta}_{j} t_i = -t_i {\zeta}_j \;\; (j \neq i,i+1),
  \\
 {\zeta}_i   {\zeta}_j
  &=&  -{\zeta}_j {\zeta}_i \;\; (i \neq j),
  \\
 {[}{\zeta}_i, e^{\eta}]
  &=& u \sum_{k \neq i}  {\text{sgn}(k-i)}
 \frac{e^\eta -e^{s_{ki} (\eta)}}{1 -e^{\text{sgn}(k-i) \cdot
 (\epsilon_k -\epsilon_i)}}[k,i].
 \end{eqnarray*}
\end{definition}
Note that the algebra $\Bn_{tr}$ is naturally a superalgebra by
declaring $e^{\pm \epsilon_i}$ to be even and ${\zeta}_i$ and
$t_i$ to be odd for every $i$. For $u \neq 0$, the subalgebra
generated by ${\zeta}_i$ $(1 \le i \le n)$ and $\Tn$ is identified
with the degenerate spin affine Hecke algebra $\Bnhat$ (cf.
Sect.~\ref{sec:affineHeckeSpin}).

\begin{theorem}   \label{th:rat=trBn}
\begin{enumerate}
 \item
 There exists a unique  homomorphism of superalgebras $\iota^-: \Bn
 \rightarrow \Bn_{tr}$, which extends the identity map on
 the subalgebra $\Tn$ such that
 \begin{eqnarray*}
 \iota^- (y_i) &=& e^{\epsilon_i} \\
 \iota^- ({\xi}_i) &=& e^{-\epsilon_i} \left( {\zeta}_i - u \sum_{k<i}
 [k,i] \right).
 \end{eqnarray*}
 Moreover, $\iota^-$ is injective.
 \item
 The $\iota^-$ extends to an algebra isomorphism:
$$\C [\underline{y}^\pm] \otimes_{\C [ \underline{y} ]} \Bn
 \stackrel{\cong}{\longrightarrow}
 \Bn_{tr}.$$
The inverse map $j^-$ is an extension of the identity map on $\Tn$
such that
 \begin{eqnarray*}
 j^- (e^{\epsilon_i}) &=& y_i,  \\
 j^- ({\zeta}_i) &=& y_i {\xi}_i + u\sum_{k<i} [k,i].
 \end{eqnarray*}
\end{enumerate}
\end{theorem}

\begin{proof}
Follows from Theorem~\ref{th:isoDAHA} and Theorem~\ref{th:rat=tr};
or one argues directly as for Theorem~\ref{th:rat=tr}.
\end{proof}

Theorem~\ref{th:rat=trBn} and the PBW Theorem~\ref{th:PBW-Bn} for
$\Bn$ immediately imply the following.
\begin{corollary}
The PBW property holds for $\Bn_{tr}$. That is, the multiplication
of the subalgebras gives rise to an isomorphism of vector spaces:
$$\C[P] \otimes \Tn \otimes \C[ {\zeta}_1,
\ldots, {\zeta}_n] \stackrel{\cong}{\longrightarrow} \Bn_{tr}.
$$
\end{corollary}
\subsection{An algebra isomorphism} The following theorem can be established by a direct computation.
\begin{theorem}
There exists a superalgebra isomorphism:
\begin{align*}
 \Phi^{tr}: & \An_{tr} \stackrel{\cong}{\longrightarrow} \Cl_n \bigotimes
 \Bn_{tr}, \\
e^{\epsilon_i} \mapsto e^{\epsilon_i}, \; c_i \mapsto c_i, & \;
 \epsilon_i^\vee \mapsto {\sqrt{-2}}c_i  {\zeta}_i, \;
 s_i \mapsto \frac{1}{\sqrt{-2}} (c_i -c_{i+1})  t_i.
\end{align*}
The inverse map is given by
\begin{eqnarray*}
e^{\epsilon_i} \mapsto e^{\epsilon_i}, \; c_i  \mapsto c_i, \;
  {\zeta}_i \mapsto \frac{1}{\sqrt{-2}} c_i \epsilon_i^\vee, \;
 t_i \mapsto  \frac{1}{\sqrt{-2}} (c_{i+1} -c_i) s_i.
\end{eqnarray*}
\end{theorem}

{\bf Notes added.} This first paper on spin Hecke algebras was
posted in arXiv: math.RT/0608074, 2006. The results of this paper
have been generalized and extended in different directions since
then in a series of papers \cite{Kh, KW1, KW2, KW3, W}.

\section*{Acknowledgment}

This research is partially supported by NSF and NSA grants. I
gratefully acknowledge the sesquicentennial Associateship from
University of Virginia which allows me to spend the semester of
Spring 2006 at MSRI, Berkeley. I thank MSRI for its stimulating
atmosphere and support, where this work was carried out. I thank
Ta Khongsap for his careful reading and corrections.

\end{document}